\documentclass[10pt,reqno]{amsart}
\usepackage{amsmath,amssymb,amsfonts}
\usepackage{mathrsfs}
\usepackage[hidelinks]{hyperref}
\theoremstyle{plain}
\newtheorem{theorem}{Theorem}[section]
\newtheorem{lemma}[theorem]{Lemma}
\newtheorem{remark}[theorem]{Remark}
\newtheorem{corollary}[theorem]{Corollary}
\numberwithin{equation}{section}
\newcommand{\R}{\mathbb R}
\newcommand{\Z}{\mathbb Z}
\newcommand{\T}{\mathbb T}
\newcommand{\Leb}{\operatorname{Leb}}
\newcommand{\dimH}{\dim_{\mathrm H}}
\newcommand{\norm}[1]{\lVert #1\rVert}
\title[DIMENSION DROP FOR UNIFORM LITTLEWOOD LEVELS]{Dimension drop for superlevel sets of the uniform Littlewood constant}
\author{Shahriar Mirzadeh}
\address{Department of Mathematics, Brandeis University,
415 South Street, Waltham, MA 02453, USA}
\email{shahmir@brandeis.edu}
\begin{document}
\begin{abstract}
For $(x,y)\in\R^2$ and $Q\ge1$, set $F_Q(x,y)=Q\min_{1\le q\le Q}\norm{qx}\norm{qy}$, where $q$ ranges over positive integers and $\norm{\cdot}$ denotes distance to the nearest integer. The set on which the uniform Littlewood constant $L=\limsup_{Q\to\infty}F_Q$ is positive has full Hausdorff dimension. We prove that for every fixed $c>0$ the superlevel set $\{L\ge c\}$ has Hausdorff dimension strictly less than two. More precisely, the Lebesgue measure of $\{(x,y)\in[0,1]^2:F_Q(x,y)\ge c\}$ is $O_c(Q^{-\kappa(c)})$ for some $\kappa(c)>0$ and all $Q\ge2$. The proof combines large deviations for one diagonal action on the space of unimodular lattices with effective equidistribution of expanding horospheres for another.
\end{abstract}
\maketitle
\section{Introduction and main results}
The classical Littlewood conjecture asserts that
\begin{equation}\label{eq:littlewood}
\liminf_{q\to\infty}q\norm{qx}\norm{qy}=0
\end{equation}
for every $(x,y)\in\R^2$, where $\norm{t}$ denotes the distance from $t$ to the nearest integer. Bandi, Fregoli and Kleinbock~[1] introduced a uniform version of this problem. For $Q\ge1$, set
\begin{equation}\label{eq:uniform}
F_Q(x,y):=Q\min_{1\le q\le Q}\norm{qx}\norm{qy},\qquad L(x,y):=\limsup_{Q\to\infty}F_Q(x,y),
\end{equation}
where the minimum is taken over positive integers $q$. The uniform Littlewood conjecture would assert that $L(x,y)=0$ for all $(x,y)\in\R^2$. Equivalently, for every $\varepsilon>0$ and all sufficiently large $Q$, one could find $1\le q\le Q$ such that
\[
\norm{qx}\norm{qy}<\frac{\varepsilon}{Q}.
\]
Unlike the classical conjecture, the uniform version requires approximations at every sufficiently large scale.

Bandi, Fregoli and Kleinbock~[1] proved that $L(x,y)=0$ for Lebesgue-almost every $(x,y)$. Schleischitz~[7] subsequently disproved the uniform conjecture and showed that its counterexample set is residual. Further constructions were given by Moshchevitin~[5] and Shulga~[8], the latter obtaining a set of counterexamples of Hausdorff dimension at least $3/2$. More recently, Neckrasov, Wu and Yang~[6] proved that
\begin{equation}\label{eq:E-intro}
E:=\{(x,y)\in\R^2:L(x,y)>0\}
\end{equation}
is hyperplane absolute winning. In particular, $\dimH E=2$. They also show that, for sufficiently rapidly growing sequences $Q_m$, the set of pairs satisfying $\liminf_{m\to\infty}F_{Q_m}(x,y)>0$ is hyperplane absolute winning~[6]. In these results, the positive threshold may depend on the pair.

In this paper we study the Hausdorff dimension of the fixed superlevel sets
\begin{equation}\label{eq:Ec}
E_c:=\{(x,y)\in\R^2:L(x,y)\ge c\},\qquad c>0.
\end{equation}
Since $E=\bigcup_{n\ge1}E_{1/n}$, the result of~[6] implies that $\sup_{c>0}\dimH E_c=2$. It does not, however, determine whether any individual set $E_c$ has full Hausdorff dimension. Moshchevitin~[5, Theorem~1] constructs a pair with
$L(x,y)\ge 2/(3+\sqrt5)$, so these sets are nonempty
for some fixed positive thresholds.
(The strict cutoff $q<Q$ in~[5] yields the same
limsup, since at integer scales the corresponding
quantity equals $\frac{Q}{Q-1}F_{Q-1}(x,y)$,
and $\frac{Q}{Q-1}\to1$.)
Our main result shows that their dimensions are
nevertheless strictly less than two.

\begin{theorem}\label{thm:main}
For every $c>0$ there exists $\delta(c)>0$ such that
\begin{equation}\label{eq:drop}
\dimH E_c\le 2-\delta(c)<2.
\end{equation}
\end{theorem}
Thus, although the set of counterexamples to the uniform Littlewood conjecture has full Hausdorff dimension, every fixed positive superlevel set has dimension strictly less than two. In particular, the full dimension of $E$ is not attained at any positive level.

Theorem~\ref{thm:main} follows from a quantitative estimate for the finite-scale sets. Since $F_Q$ is $\Z^2$-periodic, we regard it as a function on $\T^2=\R^2/\Z^2$. Let $\Leb$ denote normalized Lebesgue measure on $\T^2$, and write
\begin{equation}\label{eq:AQ}
A_Q(c):=\{w\in\T^2:F_Q(w)\ge c\}.
\end{equation}
\begin{theorem}\label{thm:measure}
For every $0<c<1$ there exist constants $C_c>0$ and $\kappa(c)\in(0,1]$ such that
\begin{equation}\label{eq:measure}
\Leb\bigl(A_Q(c)\bigr)\le C_cQ^{-\kappa(c)}\qquad(Q\ge2).
\end{equation}
Consequently,
\begin{equation}\label{eq:dim}
\dimH E_c\le 2-\frac{\kappa(c/8)}{2}.
\end{equation}
\end{theorem}
The measure estimate is nontrivial for $0<c<1/2$. Indeed, the one-dimensional Dirichlet theorem gives, for $Q\ge1$,
\begin{equation}\label{eq:dirichlet}
F_Q(x,y)\le\frac{Q}{2(\lfloor Q\rfloor+1)}<\frac12,\qquad L(x,y)\le\frac12.
\end{equation}
Thus $A_Q(c)=\varnothing$ for $c\ge1/2$, while $E_c=\varnothing$ for $c>1/2$. We state Theorem~\ref{thm:measure} for $c<1$ to avoid separate cases in the subsequent argument. The passage from Theorem~\ref{thm:measure} to Theorem~\ref{thm:main} uses a $Q^2$-Lipschitz estimate for $F_Q$ and a dyadic covering argument. The method does not provide an explicit formula for $\kappa(c)$.

Combining Theorem~\ref{thm:main} with the full-dimensionality result of~[6], we obtain the following consequence.
\begin{corollary}\label{cor:limit}
The function $c\mapsto\dimH E_c$ is nonincreasing, and
\begin{equation}\label{eq:limit}
\lim_{c\downarrow0}\dimH E_c=2.
\end{equation}
\end{corollary}
\begin{proof}
The inclusion $E_{c_2}\subset E_{c_1}$ for $c_2>c_1$ gives monotonicity. Furthermore, $E=\bigcup_{n\ge1}E_{1/n}$. Countable stability and~[6] give $\sup_n\dimH E_{1/n}=2$.
\end{proof}

\noindent\textbf{Outline of the proof.} Let
\[
X=\operatorname{SL}_3(\R)/\operatorname{SL}_3(\Z),\qquad
 a_T=\operatorname{diag}(e^{T/2},e^{T/2},e^{-T}),
\]
\[
u_w=\begin{pmatrix}1&0&x\\0&1&y\\0&0&1\end{pmatrix},\qquad
b_s=\operatorname{diag}(e^s,e^{-s},1),\quad w=(x,y).
\]
We first show that for every $0<c<1$ there is a nonempty open set $\Omega_c\subset X$ such that, for $Q=e^T$,
\begin{equation}\label{eq:avoid}
F_Q(w)\ge c\quad\Longrightarrow\quad b_sa_Tu_w\Z^3\notin\Omega_c\qquad(0\le s\le T/2).
\end{equation}
The exponential large-deviation estimate of Kim and Mor~[2],
applied to $b=b_1$, bounds the Haar measure of points whose
$b$-orbits visit $\Omega_c$ too infrequently. To estimate the
corresponding set of parameters $w\in[0,1]^2$, we must transfer
this bound to the family $\{a_Tu_w\Z^3:w\in[0,1]^2\}$.
Since the induced parameter measure is singular with respect
to Haar measure on $X$, we use effective equidistribution of
expanding horospheres~[3] to compare the two measures on
suitable smooth approximations of the exceptional event.

To combine these estimates, choose a nonnegative smooth function $\varphi$ supported in $\Omega_c$ and consider the orbit averages
\[
P_N(z)=\frac1N\sum_{j=0}^{N-1}\varphi(b^jz).
\]
We transfer the estimate by constructing a smooth majorant
of the indicator of the deviation event. Its Haar integral
is controlled by the large-deviation estimate, while its
Sobolev norm grows at most exponentially in $N$.
Effective equidistribution then gives
\begin{equation}\label{eq:transfer-outline}
\text{exceptional parameter measure}
\ll e^{-\gamma N}+e^{-\lambda T+B_0N},
\end{equation}
where $\gamma>0$ is the large-deviation rate, $\lambda>0$
is the equidistribution rate, and $B_0>0$ bounds the
exponential growth of the Sobolev norm. These constants
and the implied constant are independent of $N$ and $T$,
with the level $c$ fixed.
Choosing $N=\lfloor\eta T\rfloor$ with
$0<\eta<\min\{1/2,\lambda/B_0\}$ makes both terms decay
exponentially in $T$, yielding polynomial decay in $Q=e^T$.

The two diagonal actions serve different purposes: the $b_s$-action detects the Diophantine condition, whereas the $a_T$-action provides equidistribution of the parameter family. This differs from the single-flow dimension-drop arguments in~[3,4]. Section~2 recalls the dynamical estimates; Section~3 establishes the lattice-avoidance implication; Section~4 proves a general transfer theorem and deduces the finite-scale estimate; and Section~5 contains the Hausdorff-dimension argument.

\medskip\noindent\textbf{Acknowledgments.}
I am grateful to Dmitry Kleinbock for introducing me
to the uniform Littlewood conjecture.

\medskip\noindent\textbf{AI Declaration.}
ChatGPT 6 was used as an auxiliary tool
to check selected calculations and mathematical arguments
and to improve the exposition.
The author independently reviewed and verified the
mathematical content and takes full responsibility
for the manuscript.

\section{Dynamical preliminaries}
Let
\[
G=\operatorname{SL}_3(\R),\qquad \Gamma=\operatorname{SL}_3(\Z),\qquad X=G/\Gamma,
\]
with $G$-invariant probability measure $\mu$. We identify $g\Gamma$ with the unimodular lattice $g\Z^3$ and put $z_0=\Gamma$. For $w=(x,y)\in\R^2$, set
\begin{equation}\label{eq:actions}
u_w=\begin{pmatrix}1&0&x\\0&1&y\\0&0&1\end{pmatrix},\quad
a_T=\operatorname{diag}(e^{T/2},e^{T/2},e^{-T}),\quad
b_s=\operatorname{diag}(e^s,e^{-s},1).
\end{equation}
The actions $a_T$ and $b_s$ commute, and
\[
a_Tu_wa_T^{-1}=u_{e^{3T/2}w}.
\]
Consequently, $U=\{u_w:w\in\R^2\}$ is the full expanding
horospherical subgroup of $a_T$.

Fix a right-invariant Riemannian metric on $G$ whose induced
Riemannian metric on $U\cong\R^2$ is the standard Euclidean
metric, and equip $X$ with the quotient metric.
Normalize Haar measure on $U$ so that $du_w=dw$.

Fix a basis $V_1,\ldots,V_D$ of the Lie algebra of $G$. For a Lie-algebra vector $V$ and a smooth function $f$ on $X$, let
\[
D_V f(z)=\left.\frac{d}{dt}\right|_{t=0}f(\exp(tV)z).
\]
If $I$ ranges over compositions of these derivatives, including the identity, write
\begin{equation}\label{eq:sobolev}
\begin{aligned}
\|f\|_{C^k}&:=\sum_{|I|\le k}\|D_If\|_\infty,\qquad
\|f\|_{k,2}:=\sum_{|I|\le k}\|D_If\|_{L^2(X,\mu)},\\
S_\ell(f)&:=\max\{\|f\|_{C^1},\|f\|_{\ell,2}\}.
\end{aligned}
\end{equation}
We use the corresponding norms on other homogeneous spaces after fixing a basis of their Lie algebras.

\medskip\noindent\textbf{Large deviations for $b$.}
The group $G$ is connected and semisimple, has finite center and no compact factors, and $\Gamma$ is irreducible. The strong spectral gap gives exponential mixing for the nontrivial $\operatorname{Ad}$-diagonalizable element $b=b_1=\operatorname{diag}(e,e^{-1},1)$; see~[2, \S1.2]. Thus [2, Theorem~1.2] applies. For every real-valued $\varphi\in C_c^\infty(X)$ there is $c_\varphi>0$ such that, for $0<\varepsilon<1$ and $N\gg_\varphi\varepsilon^{-2}$,
\begin{equation}\label{eq:ld}
\mu\left\{z\in X:\left|\frac1N\sum_{j=0}^{N-1}\varphi(b^jz)-\mu(\varphi)\right|>\varepsilon\right\}
\le e^{-c_\varphi\varepsilon^2N}.
\end{equation}
Here $\mu(\varphi)=\int_X\varphi\,d\mu$. The two-sided estimate in \eqref{eq:ld} includes the lower tail, which is the one needed here: if an orbit avoids the support of a nonnegative $\varphi$ with $\mu(\varphi)>0$, then its orbit average is zero.

\medskip\noindent\textbf{Effective equidistribution along $U$.}
The flow $a_T$ is exponentially mixing and $U$ is its full expanding horospherical subgroup. It follows from [3, Definition~1.1 and Theorem~1.2] that there exist $A,B,\lambda>0$ and an integer $\ell\ge1$ such that
\begin{equation}\label{eq:eqd}
\left|\int_U f(u)\psi(a_Tuz)\,du-\left(\int_U f(u)\,du\right)\mu(\psi)\right|
\ll \|f\|_{C^\ell}S_\ell(\psi)e^{-\lambda T}
\end{equation}
for $f\in C_c^\infty(U)$ supported in the unit ball and smooth $\psi$ with finite Sobolev norms of all orders and finite $C^1$ norm, provided
\begin{equation}\label{eq:inj}
T\ge A+B\log\frac1{r_0(z)}.
\end{equation}
Here $r_0(z)$ denotes the injectivity radius of $z$. The implicit constant in \eqref{eq:eqd} does not depend on $\psi$ beyond the displayed norm. This uniformity will be used with test functions that depend on $T$.

We will apply \eqref{eq:eqd} with a parameter cutoff supported in an arbitrary fixed compact subset of $U\cong\R^2$.
\begin{lemma}\label{lem:cutoff}
For every $z\in X$ and $\rho\in C_c^\infty(\R^2)$, there are constants $C_{\rho,z},T_{\rho,z}>0$ such that
\begin{equation}\label{eq:eqd-rho}
\left|\int_{\R^2}\rho(w)\psi(a_Tu_wz)\,dw
-\left(\int_{\R^2}\rho(w)\,dw\right)\mu(\psi)\right|
\le C_{\rho,z}e^{-\lambda T}S_\ell(\psi)
\end{equation}
for $T\ge T_{\rho,z}$ and every $\psi$ admissible in \eqref{eq:eqd}. The constants are independent of $\psi$.
\end{lemma}
\begin{proof}
Choose finitely many points $w_i\in\R^2$ such that
$\operatorname{supp}\rho$ is covered by the balls
$w_i+B(0,1/2)$. Take a smooth partition
$\rho=\sum_i\rho_i$ subordinate to this cover, and set
$f_i(v)=\rho_i(w_i+v)$ and $z_i=u_{w_i}z$.
Then $f_i\in C_c^\infty(U)$ is supported in $B(0,1)$.
Since $U$ is abelian,
\[
\int_{\R^2}\rho_i(w)\psi(a_Tu_wz)\,dw
=\int_U f_i(v)\psi(a_Tu_vz_i)\,dv.
\]
Apply \eqref{eq:eqd} to each pair $(f_i,z_i)$ and sum.
The finitely many lower bounds on $T$ in
\eqref{eq:inj} are absorbed into $T_{\rho,z}$.
None of the resulting constants depends on $\psi$.
\end{proof}

\section{From Diophantine approximation to lattice avoidance}
We show that the inequality $F_Q(w)\ge c$ forces a segment of the $b_s$-orbit of $a_Tu_w\Z^3$ to avoid a fixed nonempty open set. Fix $0<c<1$, and define
\begin{equation}\label{eq:omega}
\Omega_c:=\left\{\Lambda\in X:\begin{array}{l}
\text{there is }(v_1,v_2,v_3)\in\Lambda\setminus\{0\}\text{ with}\\
|v_1|<\sqrt c,\quad |v_2|<\sqrt c,\quad |v_3|<1
\end{array}\right\}.
\end{equation}
The set $\Omega_c$ is open and nonempty: for example, $\operatorname{diag}(t,1,t^{-1})\Z^3\in\Omega_c$ whenever $0<t<\sqrt c$.
\begin{lemma}\label{lem:avoid}
If $Q=e^T$ and $w\in A_Q(c)$, then
\begin{equation}\label{eq:avoid-lemma}
b_sa_Tu_wz_0\notin\Omega_c\qquad(0\le s\le T/2).
\end{equation}
\end{lemma}
\begin{proof}
Suppose that $b_sa_Tu_w\Z^3\in\Omega_c$ for some $0\le s\le T/2$. There is $(p,r,q)\in\Z^3\setminus\{0\}$ such that
\[
\left(e^{T/2+s}(p+qx),\ e^{T/2-s}(r+qy),\ e^{-T}q\right)
\]
has coordinate absolute values less than $\sqrt c$, $\sqrt c$ and $1$, respectively. If $q=0$, then $e^{T/2+s}\ge1$ and $e^{T/2-s}\ge1$, so $|p|,|r|<\sqrt c<1$. This contradicts $(p,r,q)\ne0$. Thus $q\ne0$, and the third coordinate gives $1\le|q|<Q$. The first two coordinates give
\[
\|qx\|<\sqrt c\,e^{-T/2-s},\qquad \|qy\|<\sqrt c\,e^{-T/2+s}.
\]
It follows that $\|qx\|\|qy\|<ce^{-T}=c/Q$, contradicting $F_Q(w)\ge c$.
\end{proof}
In particular, if $w\in A_Q(c)$, then $b^ja_Tu_wz_0\notin\Omega_c$ for every integer $0\le j\le T/2$.

\section{Transfer to parameter measure}
For the application, we must estimate the set of parameters $w$
for which the $b$-orbit of $a_Tu_wz_0$ visits $\Omega_c$
infrequently. The large-deviation bound \eqref{eq:ld} controls
the corresponding set in $X$ with respect to Haar measure.
The next theorem transfers this estimate to the parameter
family using effective equidistribution. The key step is to
construct a smooth majorant of the indicator of the deviation
event and control its Sobolev norm.

We formulate the argument on a general homogeneous space.
Let $G$ be a finite-dimensional Lie group, $\Gamma<G$ a lattice,
and $X=G/\Gamma$, with $G$-invariant probability measure $\mu$.
Let $H<G$ be a closed subgroup with Haar measure $\nu$.
Fix a one-parameter semigroup $(a_T)_{T\ge0}\subset G$
and an element $b\in G$. For $x\in X$, consider the parameter
family $h\mapsto a_Thx$, where $h\in H$.
Fix a basis of the Lie algebra of $G$ and define $S_\ell$
as in \eqref{eq:sobolev}.

Assume the following two estimates.

\medskip\noindent\textbf{(E) Effective equidistribution.} There exist $\lambda>0$ and an integer $\ell\ge1$ such that for every $x\in X$, $\rho\in C_c^\infty(H)$, and $f\in C_c^\infty(X)$,
\begin{equation}\label{eq:E-hyp}
\left|\int_H\rho(h)f(a_Thx)\,d\nu(h)-\nu(\rho)\mu(f)\right|
\le C_{\rho,x}e^{-\lambda T}S_\ell(f)
\end{equation}
whenever $T\ge T_{\rho,x}$, where $\nu(\rho)=\int_H\rho\,d\nu$. The constants $C_{\rho,x}$ and $T_{\rho,x}$ are independent of $f$, including its support.

\medskip\noindent\textbf{(LD) Exponential large deviations.} For every real-valued $\varphi\in C_c^\infty(X)$ and $0<d<1$, there exist $\gamma=\gamma(\varphi,d)>0$ and $N_0=N_0(\varphi,d)$ such that
\begin{equation}\label{eq:LD}
\mu\{z\in X:|P_N^\varphi(z)-\mu(\varphi)|>d\}\le e^{-\gamma N}\qquad(N\ge N_0),
\end{equation}
where
\begin{equation}\label{eq:PN}
P_N^\varphi(z)=\frac1N\sum_{j=0}^{N-1}\varphi(b^jz).
\end{equation}
No commutation relation between $a_T$ and $b$ is assumed. The uniformity over $f$ in (E) is essential, since the compact support of the test function constructed below will depend on $N$.

\begin{theorem}[Transfer of large deviations]\label{thm:transfer}
Assume \textup{(E)} and \textup{(LD)}. There exists $\eta\in(0,1/2)$ such that for every real-valued $\varphi\in C_c^\infty(X)$ and every $\varepsilon>0$, there is $\kappa>0$ with the following property: for every $x\in X$ and compact $K\subset H$, putting $N=\lfloor\eta T\rfloor$ gives
\begin{equation}\label{eq:transfer}
\nu\{h\in K:|P_N^\varphi(a_Thx)-\mu(\varphi)|\ge\varepsilon\}
\ll_{\varphi,\varepsilon,K,x}e^{-\kappa T}
\end{equation}
for all sufficiently large $T$. The exponents $\eta$ and $\kappa$ are independent of $K$ and $x$; moreover, $\eta$ can be chosen independently of $\varphi$ and $\varepsilon$.
\end{theorem}
\begin{proof}
Fix $\varphi\in C_c^\infty(X)$ and $\varepsilon>0$, and abbreviate $P_N=P_N^\varphi$ and $m=\mu(\varphi)$. We wish to estimate the parameters $h\in K$ for which $|P_N(a_Thx)-m|\ge\varepsilon$. The indicator of this event is not smooth, so we replace it by a smooth majorant. We will then use (E) to compare its integral over the parameter family with its Haar integral.

Put $d=\min\{\varepsilon/2,1/2\}$ and choose $\chi\in C^\infty(\R,[0,1])$ such that
\[
\chi(t)=0\quad(|t|\le d),\qquad \chi(t)=1\quad(|t|\ge\varepsilon).
\]
For $N\ge1$, define $M_N(z)=\chi(P_N(z)-m)$. By construction,
\begin{equation}\label{eq:maj}
\mathbf{1}_{\{|P_N-m|\ge\varepsilon\}}\le M_N
\le\mathbf{1}_{\{|P_N-m|>d\}}.
\end{equation}
Since $d\in(0,1)$, hypothesis (LD) applied to $\varphi$ at threshold $d$ gives constants $\gamma>0$ and $N_0\ge1$ such that
\begin{equation}\label{eq:muMN}
\mu(M_N)\le e^{-\gamma N}\qquad(N\ge N_0).
\end{equation}

The function $M_N$ need not have compact support, whereas (E) is stated for compactly supported test functions. This is resolved by subtracting the constant value of $M_N$ outside a compact set. Indeed, let
\[
B_N=\bigcup_{j=0}^{N-1}b^{-j}(\operatorname{supp}\varphi).
\]
Each set in this finite union is compact. If $z\notin B_N$, then $\varphi(b^jz)=0$ for every $0\le j<N$; hence $P_N(z)=0$ and $M_N(z)=\chi(-m)$. Therefore
\begin{equation}\label{eq:fN}
c_0=\chi(-m),\qquad f_N=M_N-c_0\in C_c^\infty(X),\qquad
\operatorname{supp}f_N\subset B_N.
\end{equation}
Subtracting $c_0$ has no effect on the equidistribution error: for any $\rho\in C_c^\infty(H)$,
\begin{equation}\label{eq:cancel}
\begin{aligned}
&\int_H\rho(h)M_N(a_Thx)\,d\nu(h)-\nu(\rho)\mu(M_N)\\
&\qquad=\int_H\rho(h)f_N(a_Thx)\,d\nu(h)-\nu(\rho)\mu(f_N).
\end{aligned}
\end{equation}
We may thus apply (E) to $f_N$. Its support depends on $N$, so it remains to control the Sobolev norm by a bound whose constants do not depend on that support.

Fix a norm on the Lie algebra of $G$. Since $b$ is fixed, there is $\beta>0$ for which
\begin{equation}\label{eq:Ad}
\|\operatorname{Ad}(b^j)\|\le e^{\beta j}\qquad(j\ge0).
\end{equation}
For any Lie-algebra vector $V$, the identity $b^j\exp(tV)b^{-j}=\exp(t\operatorname{Ad}(b^j)V)$ implies
\begin{equation}\label{eq:dAd}
D_V(\varphi\circ b^j)(z)=\bigl(D_{\operatorname{Ad}(b^j)V}\varphi\bigr)(b^jz).
\end{equation}
Applying this identity successively to an ordered derivative $D_I=D_{V_{i_1}}\cdots D_{V_{i_k}}$ gives
\[
D_I(\varphi\circ b^j)=\bigl(D_{\operatorname{Ad}(b^j)V_{i_1}}\cdots D_{\operatorname{Ad}(b^j)V_{i_k}}\varphi\bigr)\circ b^j.
\]
With respect to the fixed Lie-algebra basis, every coefficient of $\operatorname{Ad}(b^j)V_{i_r}$ is $O(e^{\beta j})$. Consequently, for $k\ge1$,
\begin{equation}\label{eq:deriv}
\|D_I(\varphi\circ b^j)\|_\infty\le C_k\|\varphi\|_{C^k}e^{k\beta j},
\end{equation}
where $C_k$ depends only on the fixed basis and norm. Averaging over $0\le j<N$ yields
\begin{equation}\label{eq:derivPN}
\begin{aligned}
\|D_IP_N\|_\infty&\le\frac1N\sum_{j=0}^{N-1}\|D_I(\varphi\circ b^j)\|_\infty\\
&\le C_k\|\varphi\|_{C^k}e^{k\beta N}.
\end{aligned}
\end{equation}
For $k\ge1$, the derivatives of $f_N$ agree with those of $\chi(P_N-m)$. By the chain rule, each derivative of order $k$ is a sum of products of derivatives of $P_N$ whose total orders add to $k$, multiplied by bounded derivatives of $\chi$. Equation \eqref{eq:derivPN} therefore implies
\[
\|D_If_N\|_\infty\ll_{\varphi,\varepsilon,k}e^{k\beta N}
\qquad(1\le|I|=k\le\ell).
\]
Also $|f_N|\le1$. Since $\mu$ is a probability measure, the $L^2(\mu)$ norms of these derivatives satisfy the same bounds. Thus, with $B=\ell\beta>0$,
\begin{equation}\label{eq:SfN}
S_\ell(f_N)\le C_{\varphi,\varepsilon}e^{BN}.
\end{equation}
Here $B$ depends only on $b$, $\ell$, and the fixed choice
of norms, while $C_{\varphi,\varepsilon}$ is independent of $N$.
In particular, the varying support of $f_N$ introduces no
additional factor.

Now fix $x\in X$ and a compact set $K\subset H$. Choose a nonnegative $\rho\in C_c^\infty(H)$ with $\rho\ge1$ on $K$, and write
\[
E_{N,T}=\{h\in K:|P_N(a_Thx)-m|\ge\varepsilon\}.
\]
The first inequality in \eqref{eq:maj} and the choice of $\rho$ give
\[
\nu(E_{N,T})\le\int_H\rho(h)M_N(a_Thx)\,d\nu(h).
\]

By \eqref{eq:cancel} and (E), applied to $f_N$, we have
\[
\begin{aligned}
\int_H\rho(h)M_N(a_Thx)\,d\nu(h)
&=\nu(\rho)\mu(M_N)\\
&\quad+\int_H\rho(h)f_N(a_Thx)\,d\nu(h)-\nu(\rho)\mu(f_N)\\
&\le \nu(\rho)\mu(M_N)+C_{\rho,x}e^{-\lambda T}S_\ell(f_N).
\end{aligned}
\]
Combining this with \eqref{eq:muMN} and \eqref{eq:SfN} gives
\begin{equation}\label{eq:master}
\begin{aligned}
\nu(E_{N,T})
&\le\nu(\rho)\mu(M_N)+C_{\rho,x}e^{-\lambda T}S_\ell(f_N)\\
&\le\nu(\rho)e^{-\gamma N}+C_{\rho,x}C_{\varphi,\varepsilon}e^{-\lambda T+BN},
\end{aligned}
\end{equation}
provided $N\ge N_0$ and $T\ge T_{\rho,x}$.

To make both terms in \eqref{eq:master} decay exponentially in $T$, set
\begin{equation}\label{eq:eta}
\eta=\min\left\{\frac14,\frac{\lambda}{2B}\right\},\qquad N=\lfloor\eta T\rfloor.
\end{equation}
Then $0<\eta<1/2$ and
\[
\eta T-1\le N\le\eta T,\qquad BN\le B\eta T\le\frac{\lambda T}{2}.
\]
Thus $N\ge N_0$ for all sufficiently large $T$. Substituting these inequalities into \eqref{eq:master}, we obtain
\[
\begin{aligned}
\nu(E_{N,T})
&\le\nu(\rho)e^\gamma e^{-\gamma\eta T}
+C_{\rho,x}C_{\varphi,\varepsilon}e^{-\lambda T/2}\\
&\ll_{\varphi,\varepsilon,K,x}e^{-\kappa T},\qquad
\kappa=\min\{\gamma\eta,\lambda/2\}>0.
\end{aligned}
\]
Here $\eta$ depends only on $\lambda$ and $B$, while $\kappa$ may also depend on $\varphi$ and $\varepsilon$; neither depends on $K$ or $x$. This proves the theorem.
\end{proof}

For an open set $\Omega\subset X$, put
\[
R_N(z,\Omega)=\#\{0\le j<N:b^jz\in\Omega\}.
\]
\begin{corollary}[Low-return estimate]\label{cor:return}
Under the hypotheses of Theorem~\ref{thm:transfer}, for every nonempty open set $\Omega\subset X$, there are $\theta,\kappa>0$ such that for every $x\in X$ and compact $K\subset H$,
\begin{equation}\label{eq:return}
\nu\{h\in K:R_N(a_Thx,\Omega)\le\theta N\}
\ll_{\Omega,K,x}e^{-\kappa T},\qquad N=\lfloor\eta T\rfloor,
\end{equation}
for all sufficiently large $T$, where $\eta$ is as in Theorem~\ref{thm:transfer}. In particular, the same estimate holds for the parameters whose first $N$ iterates avoid $\Omega$.
\end{corollary}
\begin{proof}
Choose $\varphi\in C_c^\infty(X)$ with $0\le\varphi\le1$, $\varphi\not\equiv0$, and $\operatorname{supp}\varphi\subset\Omega$. Since $\mu$ has full support, $m=\mu(\varphi)>0$. Every visit contributes at most $1$ to the sum defining $P_N^\varphi$, and terms corresponding to iterates outside $\Omega$ vanish. Therefore
\[
P_N^\varphi(z)\le\frac{R_N(z,\Omega)}{N}.
\]
Take $\theta=m/4$. Whenever $R_N(z,\Omega)\le\theta N$, we have $P_N^\varphi(z)\le m/4$, so
\[
|P_N^\varphi(z)-m|\ge\frac{3m}{4}>\frac m2.
\]
The event in \eqref{eq:return} is consequently contained in the deviation event in \eqref{eq:transfer}, with $\varepsilon=m/2$. Applying Theorem~\ref{thm:transfer} proves the claimed estimate. An orbit that avoids $\Omega$ has $R_N=0$, so the last assertion follows as well.
\end{proof}

\medskip\noindent\textbf{Proof of the measure estimate in Theorem~\ref{thm:measure}.}
Identify $H=U$ with $\R^2$ through $w\mapsto u_w$; the Haar measure $\nu$ then corresponds to $dw$. Lemma~\ref{lem:cutoff} verifies (E), while \eqref{eq:ld} verifies (LD) for $b=b_1$. Fix $0<c<1$. We apply Corollary~\ref{cor:return} with $\Omega=\Omega_c$, $x=z_0$, and the compact parameter set
\[
K=\{u_w:w\in[0,1]^2\}.
\]
Let $Q=e^T$ and $N=\lfloor\eta T\rfloor$. Because $0<\eta<1/2$, we have $1\le N\le T/2$ when $T$ is sufficiently large. If $w\in A_Q(c)$, Lemma~\ref{lem:avoid} implies that
\[
b^ja_Tu_wz_0\notin\Omega_c\qquad(0\le j<N).
\]
Thus $R_N(a_Tu_wz_0,\Omega_c)=0$, and the exceptional parameters in $A_Q(c)$ form a subset of the no-return event controlled by Corollary~\ref{cor:return}. Consequently,
\begin{equation}\label{eq:measure-final}
\begin{aligned}
\Leb(A_Q(c))
&\le\nu\{u_w\in K:R_N(a_Tu_wz_0,\Omega_c)\le\theta N\}\\
&\ll_c e^{-\kappa T}=Q^{-\kappa}.
\end{aligned}
\end{equation}
The first line uses the identification of $[0,1]^2$ with a fundamental domain for $\T^2$; its boundary has measure zero. Enlarging the constant accommodates the bounded range $2\le Q<e^{T_0}$, where $T_0$ is the threshold above which the argument applies. Replacing $\kappa$ by $\min\{\kappa,1\}$ proves \eqref{eq:measure}. The dimension assertion follows from Lemma~\ref{lem:dimension} below.\hfill$\square$

\section{Hausdorff dimension bounds}
On $\T=\R/\Z$, let $d_\T(x,x')=\min_{k\in\Z}|x-x'-k|$, and equip $\T^2$ with
\[
d_\infty\bigl((x,y),(x',y')\bigr)=\max\{d_\T(x,x'),d_\T(y,y')\}.
\]
Write $\widetilde E_c=\{w\in\T^2:L(w)\ge c\}$.

\begin{lemma}\label{lem:lipschitz}
For $Q\ge1$ and $w,w'\in\T^2$,
\begin{equation}\label{eq:lipschitz}
|F_Q(w)-F_Q(w')|\le Q^2 d_\infty(w,w').
\end{equation}
\end{lemma}
\begin{proof}
For $1\le q\le Q$, the function $x\mapsto\|qx\|$ is $q$-Lipschitz on $\T$ and takes values in $[0,1/2]$. Therefore
\[
\begin{aligned}
&|\|qx\|\|qy\|-\|qx'\|\|qy'\||\\
&\qquad\le\frac q2\bigl(d_\T(x,x')+d_\T(y,y')\bigr)\le qd_\infty(w,w').
\end{aligned}
\]
The minimum of finitely many functions with the same Lipschitz constant is again Lipschitz with that constant. Hence $w\mapsto\min_{1\le q\le Q}\|qx\|\|qy\|$ is $Q$-Lipschitz, and multiplication by $Q$ gives \eqref{eq:lipschitz}.
\end{proof}

\begin{lemma}\label{lem:dimension}
Suppose that for every $0<a<1$ there exists $\kappa(a)\in(0,1]$ such that
\begin{equation}\label{eq:measure-hyp}
\Leb(A_Q(a))\ll_a Q^{-\kappa(a)}\qquad(Q\ge2).
\end{equation}
Then, for every $0<c<1$,
\[
\dimH E_c\le2-\frac{\kappa(c/8)}2.
\]
\end{lemma}
\begin{proof}
We first convert the measure estimate into a covering estimate for $A_Q(a)$ at spatial scale $Q^{-2}$. Fix $0<a<1$ and $Q\ge2$, and partition the unit square into $M^2$ congruent squares of side length $h=1/M$, where
\[
M=\left\lceil\frac{4Q^2}{a}\right\rceil.
\]
This partition induces a square grid on $\T^2$. If a square $S$ meets $A_Q(a)$, choose $w\in S\cap A_Q(a)$. For any $w'\in S$, we have $d_\infty(w,w')\le h$, and Lemma~\ref{lem:lipschitz} gives
\[
F_Q(w')\ge F_Q(w)-Q^2h\ge a-\frac a4\ge\frac a2.
\]
Thus the whole square $S$ lies in $A_Q(a/2)$. The interiors of the grid squares are disjoint, so if $\mathcal N_Q(a)$ denotes the number of squares meeting $A_Q(a)$, then
\[
\mathcal N_Q(a)h^2\le\Leb(A_Q(a/2)).
\]
Since $h^{-1}=M\asymp_a Q^2$, hypothesis \eqref{eq:measure-hyp} implies
\begin{equation}\label{eq:cover}
\mathcal N_Q(a)\ll_a Q^{4-\kappa(a/2)}.
\end{equation}
We have therefore covered $A_Q(a)$ by at most $O_a(Q^{4-\kappa(a/2)})$ squares, each of diameter $O_a(Q^{-2})$.

We now apply this estimate at dyadic scales. Let $w\in\widetilde E_c$. By definition of $L(w)$, there are arbitrarily large $Q$ with $F_Q(w)\ge c/2$. Choose the integer $n$ such that $2^n\le Q<2^{n+1}$. The set of denominators $1\le q\le2^n$ is contained in the set $1\le q\le Q$, so
\[
\begin{aligned}
F_{2^n}(w)
&=2^n\min_{1\le q\le2^n}\|qx\|\|qy\|\\
&\ge\frac{2^n}{Q}F_Q(w)\ge\frac c4.
\end{aligned}
\]
As such values of $Q$ can be arbitrarily large, so can the corresponding $n$. It follows that
\begin{equation}\label{eq:limsup}
\widetilde E_c\subset\limsup_{n\to\infty} A_{2^n}(c/4).
\end{equation}

Apply \eqref{eq:cover} with $a=c/4$ and $Q=2^n$. For each $n$, the set $A_{2^n}(c/4)$ can be covered by at most
\[
O_c\bigl(2^{n(4-\kappa(c/8))}\bigr)
\]
squares of diameter $O_c(2^{-2n})$. For any $s>2-\kappa(c/8)/2$, the sum of the $s$th powers of these diameters over all $n\ge n_0$ is bounded by a constant times
\[
\sum_{n\ge n_0}2^{n(4-\kappa(c/8)-2s)}.
\]
Its exponent is negative, so this tail tends to zero as $n_0\to\infty$. In view of \eqref{eq:limsup}, the definition of Hausdorff dimension yields
\[
\dimH\widetilde E_c\le2-\frac{\kappa(c/8)}2.
\]
Finally, $E_c$ is $\Z^2$-periodic and is the full inverse image of $\widetilde E_c$ under the quotient map $\R^2\to\T^2$. On sufficiently small coordinate squares this map is bi-Lipschitz, and both spaces can be covered by countably many such squares. Since Hausdorff dimension is preserved under bi-Lipschitz maps and countable unions take the supremum of the component dimensions,
$\dimH E_c=\dimH\widetilde E_c$, proving the lemma.
\end{proof}

\medskip\noindent\textbf{Completion of
Theorems~\ref{thm:main} and~\ref{thm:measure}.}
The measure bound in Theorem~\ref{thm:measure}
and Lemma~\ref{lem:dimension} give \eqref{eq:dim},
and hence Theorem~\ref{thm:main}, for $0<c<1$.
For $c>1/2$, we have $E_c=\varnothing$ by
\eqref{eq:dirichlet}. This completes the proofs.
\hfill$\square$

\begin{remark}[Products in higher dimensions]
For $d\ge2$, put
\[
L_d(x_1,\ldots,x_d)=\limsup_{Q\to\infty}Q\min_{1\le q\le Q}\prod_{i=1}^d\|qx_i\|,
\qquad E_c^{(d)}=\{x\in\R^d:L_d(x)\ge c\}.
\]
Since $\|qx_i\|\le1/2$,
\[
L_d(x)\le2^{-(d-2)}L(x_1,x_2),\qquad
E_c^{(d)}\subset E_{2^{d-2}c}\times\R^{d-2}.
\]
Theorem~\ref{thm:main} therefore gives $\dimH E_c^{(d)}\le d-\delta(2^{d-2}c)<d$. This is an immediate consequence of the two-coordinate result; it does not address the dimension of $\{L_d>0\}$ or approximation by matrices.
\end{remark}

\begin{remark}[Dependence on the level]
Let $\Delta(c)=2-\dimH E_c$. By Theorem~\ref{thm:main} and Corollary~\ref{cor:limit}, $\Delta(c)>0$ for all $c>0$ and $\Delta(c)\to0$ as $c\downarrow0$. For $0<c<1$, our argument gives the implicit bound $\Delta(c)\ge\kappa(c/8)/2$. An explicit rate would require quantitative control of the large-deviation constants for smooth functions supported in the shrinking targets $\Omega_c$.
\end{remark}

\end{document}